\documentclass[preprint,12pt]{elsarticle}

\usepackage{graphicx}
\usepackage{prodint}
\usepackage{amssymb}
\usepackage{amsthm}
\usepackage{amsmath}
\usepackage{verbatim}
\usepackage{color}
\usepackage{subcaption}

\journal{Elsevier}

\begin{document}

\begin{frontmatter}



\title{A new non-parametric estimator of the cumulative distribution function under time- and random-censoring}


\author[inst1]{N. Balakrishnan}
\affiliation[inst1]{organization={Department of Mathematics and Statistics, McMaster University},
            city={Hamilton, ON},
            country={Canada, bala@mcmaster.ca} }

\author[inst2]{Ch. Paroissin}
\author[inst3]{M. Pereda Vivo}
\affiliation[inst2]{organization={Universite de Pau et des Pays de l’Adour, E2S UPPA, CNRS, LMAP},
            country={Pau, France, christian.paroissin@univ-pau.fr}}
        \affiliation[inst3]{organization={Universite de Pau et des Pays de l’Adour, E2S UPPA, CNRS, LMAP},
        	country={Pau, France, mpvivo@univ-pau.fr}}

\begin{abstract}
In this paper, we first provide a review of different non-parametric estimators for the cumulative distribution function under left-censoring. We then propose a new estimator based on a non-parametric likelihood approach using reversed hazard rate. Finally, we conclude with an application to a real data.
\end{abstract}
		
\begin{keyword}
Left-censoring \sep  Limit of dection (LOD) \sep Non-parametric likelihood method \sep  Reversed hazard rate
\MSC 62G05 \sep 62N01
\end{keyword}

\end{frontmatter}


{\em This project has received funding from the European Union’s Horizon 2020 research and innovation programme under the Marie Skłodowska-Curie grant agreement No 945416.}

	\section{Introduction}

When dealing with data analysis, we often have to deal with some censoring situations which arise when, for some units, one has only partial information. For instance, when dealing with lifetime data, some duration may not be observed exactly since the event occurs later than a certain time point. A typical case is when one perform a medical study over a given period: all the lifetimes longer than this period then get censored. Such a situation is known as right-censoring and it has been investigated rather extensively in the literature. Sometimes, censoring may occur on the left. For instance, when dealing with concentration measurements with an analytical method, one will observe an exact measurement only if it is larger than a certain threshold, called limit of detection; otherwise, one has only the information that the concentration lies between zero and this limit. Such a situation is called left-censoring. However, statistical methods and models for data subject to left censoring have received comparatively less attention in the literature. 

In this paper, we consider both time- and random left-censoring schemes. For $n$ experimental units, the quantity of interest (lifetime, concentration, etc.) is observed exactly only if it is greater than a certain threshold value. It is assumed that the observations are independent and are drawn from the same unknown distribution. Let $T_1, \ldots, T_n$ be a sample from an unknown underlying distribution $F$. Let $C_1, \ldots, C_n$ be the censoring values: these could be either deterministic in the case of time left-censoring or random in the case of random left-censoring. In the former case, the censors may be all equal or may not be equal (for instance, if there is multiple sources of censoring). In the latter case, $C_1, \ldots, C_n$ are assumed to be sample from an unknown underlying distribution, and are independent of $T_1, \ldots, T_n$. Thus, observations are $(X_1, \Delta_1), \ldots, (X_n, \Delta_n)$, where 
$$
\forall i \in \{1, \ldots, n\}, \quad
X_i = \max(T_i, C_i)
\quad {\mbox{and}} \quad
\Delta_i = \mathbb{I}_{T_i \geqslant C_i}.
$$ 

To deal with left-censored data, an easy, but a naive, approach involves replacing the censored data by anything between $0$ and this censored value. To fix the idea, let us consider the case of measuring concentrations. The used instrumentation may not provide exact value if it is below a certain known level, called limit of detection (LOD). Assume that there is a single LOD. Then, the main substitution methods are the following ones: replace any observation below the LOD by $0$, by LOD/2, by LOD/$\sqrt{2}$ or by LOD (see, for instance, \citet{HornungReed}). Of course, if one wishes to estimate the mean concentration, the first method will clearly under-estimate it, while the last method will over-estimate it. In 2010, \citet{Helsel10} recommended against using such an approach.

Since most of the papers deal with right-censoring, some alternative solutions have been proposed in the literature by transforming the data in order to switch left- to right-censoring. For instance, the following transformations have been considered: (a) $Y_i = A - X_i$ with $A$ large enough; (b) $Y_i = 1/X_i$; (c) $Y_i = - X_i$. Here, we prefer to consider a direct analysis of data with left-censoring. For the right-censoring case, the most popular non-parametric estimator of the survival function has been the one due to Kaplan and Meier \cite{KM}. Mimicking the construction of this estimator, several authors have proposed the so-called product-limit estimator for the cumulative distribution function under left-censoring situations. As we will see, some of these papers contain some mistakes. Further, some researchers  have derived a non-parametric estimator for the cumulative distribution function by using a counting process approach. 

The rest of this paper proceeds as follows. In Section 2, we introduce several notations that will be used subsequently. Section 3 is devoted to a review of some  non-parametric estimators for the cumulative distribution function under time- and random left-censoring. A pointwise estimation of the variance of the estimator is also provided. In Section 4, we introduce a new non-parametric estimator for the cumulative distribution function based on non-parametric likelihood function. This estimator is then compared to the existing ones. Finally, in Section 6, a real-life data is analysed using the proposed estimators.

\section{Notations}

In this section, we introduce some notations that we will be used in the sequel.
\begin{itemize}
	\item $n$ is the number of observations (exact or left-censored);
	\item $T_1, \ldots, T_n$ are the exact measurements (not always observed); 
	\item $C_1, \ldots, C_n$ are the censoring values (not always observed);
	\item $X_1, \ldots, X_n$ are the observed values (exact or left-censored);
	\item $\delta_1, \ldots, \delta_n$ are the indicators of the observation of exact values;
	\item $m$ is the number of distinct (exact or not) observations;
	\item $x_{(1)} < \cdots < x_{(n)}$ are the ordered distinct (exact or not) observations:
	$$
	x_{(1)} = \min\{x_i\} < x_{(2) } < \cdots < x_{(m)} = \max\{x_i\};
	$$
	\item for any $k \in \{1, \ldots, m\}$, $d_k$ is the number of exact and observed measurements equal to $x_{(k)}$:
	$$
	d_k = \# \{ i \in \{ 1, \ldots, n\} : x_i = x_{(k)} \; {\mbox{and}} \; \delta_i = 1 \};
	$$
	\item for any $k \in \{1, \ldots, m\}$, $q_k$ is the number of left-censored and observed measurements equal to $x_{(k)}$:
	$$
	q_k = \# \{ i \in \{ 1, \ldots, n\} : x_i = x_{(k)} \; {\mbox{and}} \; \delta_i = 0 \};
	$$
	\item for any $k \in \{1, \ldots, m\}$, $y_k$ is the number of observations less than or equal to $x_{(k)}$:
	$$
	y_k = \# \{ i \in \{ 1, \ldots, n\} : x_i \leqslant x_{(k)}  \} = \sum_{j=1}^k (d_j + q_j);
	$$
	\item $l$ is the number of distinct exact observations;
	\item $x_{(1)}^\star < \cdots < x_{(l)}^\star$ are the ordered distinct exact observations: 
	$$
	x_{(1)}^\star = \min\{x_i \; ; \; \delta_i = 1\} < x_{(2)}^\star < \cdots < x_{(l)}^\star = \max\{x_i \; ; \;  \delta_i = 1\};
	$$
	\item for any $k \in \{1, \ldots, l\}$, $d_k^\star$ is the number of exact measures equal to $x_{(k)}^\star$:
	$$
	d_k^\star = \# \{ i \in \{ 1, \ldots, n\} : x_i = x_{(k)}^\star \; {\mbox{and}} \; \delta_i = 1 \};
	$$
	note that $d_k^\star>0$ since, by definition, there is at least one exact observation equal to $x_{(k)}^\star$; 
	\item for any $k \in \{1, \ldots, l\}$, $y_k^\star$ is the number of observations less than or equal to $x_{(k)}$:
	$$
	y_k^\star = \# \{ i \in \{ 1, \ldots, n\} : x_i \leqslant x_{(k)}^\star  \}; 
	$$
	\item for any $k \in \{1, \ldots, l\}$, $\tilde{d}_k^\star$ is the number of (exact or not) measures equal to $x_{(k)}^\star$:
	$$
	\tilde{d}_k^\star = \# \{ i \in \{ 1, \ldots, n\} : x_i = x_{(k)}^\star \; {\mbox{and}} \; \delta_i = 1 \} \geqslant d_k^\star.
	$$
\end{itemize}

\section{Review of some non-parametric estimators of cumulative distribution function}

For time- and random-censored data, many different non-parametric estimators have been discussed for the cumulative distribution function (CDF) in the literature.  Here, we present a brief review of estimators detailing two different approaches,  with the first one being based on the chain rule and the second being based on counting processes. As we will see, these two approaches  lead to the same estimator.

\subsection{Estimator(s) based on the chain rule}

Using chain rule, we can express the cumulative distribution function at point $x_{(j)}$ with $j \in \{1, \ldots, m-1\}$, based on points $x_{(j+1)}, \ldots, x_{(m)}$, as follows:
$$
F(x_{(j)}) = \mathbb{P}[ T \leqslant x_{(j)}]  = \prod_{k=j}^{m-1} \mathbb{P}[ T \leqslant x_{(k)} | T \leqslant x_{(k+1)}] = \prod_{k=j}^{m-1} p_k.
$$
Hence, a natural estimator can be obtained by replacing $p_1, \ldots, p_{m-1}$ by some estimators $\widehat{p}_1, \ldots, \widehat{p}_{m-1}$ (since $p_1 + \cdots + p_m = 1$, $\widehat{p}_m = 1 - \widehat{p}_1 - \cdots - \widehat{p}_{m-1}$). Following the seminal work of Miller \cite{Miller}, Blackwood \cite{Blackwood} proposed the following estimator for $p_j$:
$$
\widehat{p}_j = \left(1-\frac{d_j}{y_j}\right)^{\delta _{(j)}},
$$
where $y_j$ and $d_j$ are as defined in Section 2, and  
$$
\delta_{(j)} 
= \left\{
\begin{array}{ll}
	1 &   \quad \text{if at least one observation at time } x_{(j)} \text{ is uncensored} \\
	0 &  \quad \text{if all observations at time } x_{(j)} \text{ are censored}.
\end{array}
\right.
$$
Notice that $\delta_{(j)} = 1$ (resp. $\delta_{(j)} = 0$) if, and only if, $d_j \geqslant 1$ (resp. $d_j = 0$). Blackwood \cite{Blackwood} claimed that $\widehat{p}_j$ is the maximum likelihood estimator of $p_j$ (since it is the case for the right-censoring situation). Following what is classically done for the right-censoring case (see, for instance, \cite{Miller}), we can consider $d_j$ to be a realization of a random variable $D_j$, assumed to be binomially distributed with parameters $y_j$ (observed) and $p_j$ (unobserved and unknown). It leads to the following estimator of the CDF:
\begin{equation}\label{eqn:CDF.BL}
	\forall t \geqslant 0, \quad 
	\widehat{F}(t)= \prod_{j ; x_{(j)}>t} \left(1-\frac{d_j}{y_j}\right)^{\delta _{(j)}}.
\end{equation}
As $\delta_{(j)} = 0$ is equivalent to $d_j = 0$, then this estimator can be expressed only by considering unique uncensored observations. Using the notations introduced earlier in Section 2, we have
\begin{equation}\label{eqn:CDF.BL.bis}
	\forall t \geqslant 0, \quad 
	\widehat{F}(t)= \prod_{k ; x_{(k)}^\star>t} \left(1-\frac{d_k^\star}{y_k^\star}\right).
\end{equation}

Blackwood \cite{Blackwood} also presented an estimator of the variance by applying the same approach as the one leading to the Greenwood formula for the right-censoring case. He provided the following expression:
\begin{equation}\label{eqn:VAR.BL}
	\widehat{var} \left[ \widehat{F}(t) \right] =  \left[  \widehat{F}(t) \right]^2 \sum_{j; x_{(j)}>t} \frac{d_j \delta_{(j)}}{y_j(y_j-d_j)}.
\end{equation}
It should be mentioned that Blackwood \cite{Blackwood} also proposed non-parametric estimators for the means and the quantiles based on the non-parametric estimator of the CDF. 

Despite this work of Blackwood \cite{Blackwood}, some researchers from other fields such as environmental science and physics have come up with the same estimator. Besides, Pajek {\em et al.} \cite{Pajek04a} have also considered another estimator based on a non-parametric estimator of the cumulative hazard function (CHF), such as the Nelson-Aalen and the Harrington-Flemming estimators. Unfortunately, this work contains some mistakes. It starts with an imprecise definition of the cumulative hazard function $\Lambda$ (bounds of the integral are not given and a confusion caused by using the same variable for the function $\Lambda$ and the integrand). According to Equation~(5) in \cite{Pajek04a} providing a non-parametric estimator $\widehat\Lambda$  of $\Lambda$, it seems that, in fact, they are rather considering the cumulative reversed hazard function (CRHF). It can be seen through the fact that $\widehat\Lambda$ is a decreasing function, which is not the case for the CHF, but is true in fact for the CRHF. As a consequence, the estimator proposed in Equation~(6)  in \cite{Pajek04a} is incorrect. Indeed, they have used the relationship between the CHF and the survival function (as is usually done in the right censoring case, from Nelson-Aalen estimator to Harrington-Flemming estimator). But, since they have in fact an estimator of the CRHF, they should have used the relationship between the CRHF and the CDF. This way, the correct estimator should be (using our notations) as
$$
\forall t \geqslant 0, \quad 
\widehat{F}(t) =  \prod_{k ; x_{(k)}^\star>t} \exp  \left(- \frac{d_k^\star}{y_k^\star}\right).
$$

\subsection{Estimator based on counting processes}

To the best of our knowledge, it seems that Gomez {\em et al.} \cite{Gomez92} (see also \cite{Gomez94}) were the first to propose an estimator based on counting processes. Quite surprisingly, they have derived a non-parametric estimator for the survival function of $T$ (and not for the cumulative distribution function which is more natural to consider when dealing with left censored data). They have expressed the survival function as an integral equation. By replacing other unknown functions involved in this integral equation by their empirical counterparts, they defined an estimator of the survival function of $T$. In this way, this estimator appears to be the solution of a backward Doléans equation for which the solution can be determined explicitly in the present case. Later on, Tressou \cite{Tressou} proposed a new formulation of this estimator using the formalism developed by Gill and Johansen \cite{GillJohansen90}. Tressou \cite{Tressou} considered only the case of random left-censoring and denotes by $G$ the CDF for the censoring part. For any $t \geqslant 0$, let 
$$
\mathbb{H}_n = \frac{1}{n} \sum_{i=1}^n \mathbb{I}_{x_i \leqslant t}
\quad {\mbox{and}} \quad
\mathbb{H}_{1,n} = \frac{1}{n} \sum_{i=1}^n \mathbb{I}_{x_i \leqslant t ; \delta_i = 1}
$$
be, respectively, the empirical versions of $H(t) = \mathbb{P}[X \leqslant t]$, the CDF of the (uncensored or censored) observations, and $H_1(t) = \mathbb{P}[X \leqslant t; \Delta=1]$, the CDF of the uncensored observations. The reversed hazard rate can be defined as
$$
R(t) = \int_{]t , \infty]} \frac{\mathrm{d}F}{F} =  \int_{]t , \infty]} \frac{\mathrm{d}H_1}{H}.
$$
Now, using the product integral function $\Psi$, we have $F = \Psi(R)$; see \cite{GillJohansen90}. It then follows that a non-parametric estimator of $F$ is given by
$$
\widehat{F} 
= \Prodi_{]\cdot,\infty]} \left( 1 - \mathrm{d}\widehat{R} \right)
= \Prodi_{]\cdot,\infty]} \left( 1 - \frac{\mathrm{d}\mathbb{H}_{1,n}}{\mathbb{H}_n} \right).
$$
From the above expression for $\mathbb{H}_n$ and for $\mathbb{H}_{1,n}$, we get
$$
\forall t \geqslant 0, \quad
\widehat{F}(t) = \prod_{k=1}^{l} \left( 1 - \frac{d_k^\star}{y_k^\star} \right)^{\mathbb{I}_{x_{(k)}^\star > t}}.
$$
Notice that this estimator coincides with the one given by Patilea and Rolin \cite{PatileaRolin} under the double censoring scheme when there are no right-censored observations. Using the framework of Gill and Johansen \cite{GillJohansen90}, we can derive an estimate for the variance of $\widehat{F}(t)$ as
$$
\widehat{{\mbox{var}}}\left( \widehat{F}(t) \right) = \left[ \widehat{F}(t) \right]^2 \sum_{k=1}^{l} 
\frac{d_k^\star \mathbb{I}_{x_{(k)}^\star > t}}{y_k^\star (y_k^\star - d_k^\star)}.
$$
Observe that this estimator based on counting processes is indeed the same as the one obtained by using the chain rule. 

\section{A new non-parametric estimator}

Our goal here is to develop an estimator of the CDF based on non-parametric likelihood function; it has been done for the survival function in the case of right-censoring, but not an estimator for the CDF based on left-censoring. After providing an expression of the non-parametric likelihood function in terms of reversed hazard rate (RHR), we derive an estimator for the CDF. 

Let ${\mbox{data}} = \{ (x_1, \delta_1), \ldots, (x_1, \delta_n) \}$ be the set of all observations. We can then write the likelihood function as 
$$
L(p_1, \ldots, p_m ; {\mbox{data}}) = \prod_{i=1}^{n} p_{x_i}^{\delta_i} F_{x_i}^{1-\delta_i}, 
$$
where $F_k = p_1 + \cdots + p_k$ is the CDF of the discrete distribution. Let us recall that $F_0 = 0$ and $F_m = 1$. Using the notations introduced in Section~2, we can now express $L$ as 
\begin{equation}\label{eqn:lik1}
	L(p_1, \ldots, p_m ; {\mbox{data}}) 
	= \prod_{k=1}^{m} p_k^{d_k} F_k^{q_k} 
	= \prod_{k=1}^{m} p_k^{d_k} \left( \sum_{j=1}^k p_j \right)^{q_k}.
\end{equation}
However, this expression is not tractable for optimizing with respect to $p_1, \ldots, p_m$. Instead of expressing the likelihood function in term of mass probabilities, we will rather use the notion of RHR defined as in \cite{AER}:
$$
\forall k \in \{1, \ldots, m\}, \quad 
r_k = \frac{\mathbb{P}[T = x_{(k)}]}{\mathbb{P}[T \leqslant x_{(k)}]} = \frac{p_k}{F_k}.
$$ 
Note that we have $r_1=1$ and $r_m = p_m$. As we have $p_k = F_k - F_{k-1}$, with the convention that $F_0 = 0$, we have the following relationship:
$$
\forall k \in \{1, \ldots, m\}, \quad 
r_k = \frac{F_k - F_{k-1}}{F_k} = 1 - \frac{F_{k-1}}{F_k}.
$$
By induction, we then obtain 
$$
F_{k-1} 
=  (1-r_k) F_k 
=  (1-r_k) (1-r_{k+1}) F_{k+1} 
= \cdots 
=  \prod_{j=k}^{m} (1-r_j)
$$
since $F_m=1$. We can now rewrite the likelihood function with respect to $r_2, \ldots, r_m$ (keeping in mind that $r_1=1$). Because $p_k = r_k F_k$, Equation~(\ref{eqn:lik1}) turns to be
$$
L(r_2, \ldots, r_m ; {\mbox{data}}) 
	= \prod_{k=2}^{m} r_k^{d_k} F_k^{d_k + q_k} 
	= \prod_{k=2}^{m-1} r_k^{d_k} \left[ \prod_{j=k+1}^{m} (1-r_j)\right] ^{d_k + q_k} 
	= \prod_{k=2}^{m} r_k^{d_k} (1-r_k)^{y_{k-1}}.
$$
Hence, the log-likelihood function takes os the following form:
\begin{equation*}
	\ell(r_2, \ldots,r_m; {\mbox{data}}) = \sum_{k=2}^m \{d_k\log r_k + y_{k-1} \log(1-r_k)\}.
\end{equation*}
It turns that the value of $r_k$ that maximizes the log-likelihood function to be
$$
\forall k \in \{2, \ldots, m\}, \quad 
\widehat{r}_k = \frac{d_k}{d_k + y_{k-1}} = \frac{d_k}{y_k-q_k}
$$
since $y_k = y_{k-1} + d_k + q_k$. We thus obtain the following estimator for $F_{k-1}$:
$$
\widehat{F}_{k-1} 
= \prod_{j=k}^m  \left(1-\frac{d_j}{d_j + y_{j-1}} \right) 
= \prod_{j=k}^m  \frac{y_{j-1}}{d_j + y_{j-1}} 
= \prod_{j=k}^m  \left(1-\frac{d_j}{y_j-q_j}\right),
$$ 
using which we obtain the following non-parametric estimator of the CDF:
\begin{equation}\label{estim}
	\forall t \geqslant 0, \quad \widehat{F}^{(1)}(t) = 
	\prod_{\substack{j ; x_{(j)}>t} }  \left(1-\frac{d_j}{d_j + y_{j-1}}\right) 
\end{equation}
with the convention that $\prod_\emptyset = 0$ (meaning that if $t <x_{(1)}$, then $\widehat{F}^{(1)}(t) = 0$). Observe that for all $j$ such that $d_j = 0$ (and, of course, $q_j>0$), then $1 - d_j/(d_j+y_{j-1}) = 1$. Thence, the estimator can be defined only at points $x_{(1)}^\star, \ldots, x_{(l)}^\star$, and we then obtain
\begin{equation}\label{estim_bis}
	\forall t \geqslant 0, \quad \widehat{F}^{(1)}(t) 
	= \prod_{\substack{k ; x_{(k)}^\star>t} }  \left(1-\frac{d_k^\star}{d_k^\star + y_{k-1}^\star}\right)
	= \prod_{\substack{k ; x_{(k)}^\star>t} }  \left(1-\frac{d_k^\star}{y_k^\star - q_k^\star}\right).
\end{equation}
The last expression can be interpreted as follows. As claimed in some papers cited in Section~2 (see \cite{Gillespie} and \cite{Popovic}, for instance), if there is a tie between a censored and an uncensored observations, then it is assumed that the censored value is slightly smaller than the uncensored value. In such a case, when considering $y_k^\star - q_k^\star$, we remove these censored observations of the set of individuals at-risk.

Let us compare this new estimator with the one reviewed in the last section (recall that the two approaches considered previously lead to the same estimator). Let $t \geqslant 0$ be fixed and let us then consider the ratio  
$$
\frac{\widehat{F}^{(1)}(t)}{\widehat{F}(t)}
= \prod_{\substack{k ; x_{(k)}^\star>t} }  \left(1-\frac{d_k^\star}{d_k^\star + y_{k-1}^\star}\right) /  \left(1-\frac{d_k^\star}{y_k^\star}\right) 
$$
Because $y_{k-1}^\star + d_k^\star = y_k^\star - q_k^\star \leqslant y_k^\star$ for any $ k \in \{1, \ldots, l\}$, one can easily see that $\widehat{F}^{(1)}(t) \leqslant \widehat{F}(t)$. Let us consider the special case when there are no censored measurements. In such a case, we have $q_k^\star = 0$ for all $k \in \{ 1, \ldots, l\}$ (and $l = m$). In this case, the two estimators, $\widehat{F}^{(1)}$ and $\widehat{F}$ are identical. 

The factors in the two products defining the former estimator and this new estimator differ only at points where there is both censored and uncensored measurements. It can be expected that this may occur essentially when dealing with random censoring (and with rounded values). Because products are defined from right to left, these two estimators will be different only on the lower tail. However, for left-censoring, the main issue is to estimate accurately the left tail. 

We now seek an estimator for the variance of $\widehat{F}_T^{(1)}(t)$ for a given value of $t$. For this, we assume that $(\widehat{r}_2, \ldots, \widehat{r}_m)$ is an asymptotically normal estimator of $(r_2,\ldots,r_m)$, with asymptotic covariance matrix equal to the inverse of the Fisher information. For every $k\in\{2, \ldots, m\}$, we have
\begin{equation*}
	\frac{\partial^2 \ell}{\partial r_k ^2}(\widehat{r}_2,\ldots,\widehat{r}_m;{\mbox{data}})= -\frac{(y_k-q_k)^3}{d_k(y_k-d_k-q_k)}.
\end{equation*}
So, we can conclude that 
\begin{equation*}
	\mbox{var}[\widehat{r}_k] \approx \frac{d_k(y_k-d_k-q_k)}{(y_k-q_k)^3}.
\end{equation*}
Using classical approximations for the variance based on the delta method, one can get that 
\begin{equation*}
	\widehat{{\mbox{var}}}\left[ \widehat{F}^{(1)}(t) \right] =   \left[ \widehat{F}^{(1)}(t)\right]^2 \sum_{k: X_{(k)>t}} \frac{d_k}{y_{k-1}(y_k-q_k)}.
\end{equation*}


Note that, as for the Nelson-Aalen estimator of the CHF, one can use the above results to derive a non-parametric estimator of the CRHF and deduce another non-parametric estimator of the CDF (corresponding to Harrington-Flemming estimator in the right-censoring case) thanks to the relation between CRHF and CDF.

\section{Application to a real-life data}

In this section, we use the different estimators discussed in the previous sections to analyze a real data relating to pollutants in water, with measurements being subject to left-censoring with one or multiple limit of detection (LOD) values. Here, we consider copper concentrations in shallow groundwater samples from a Basin-Through zone in the San Joaquin Valley, California (see \cite{Millard}), while studying groundwater quality. This dataset includes five different limits of detection: 1, 2, 5, 10 and 15. There are multiple limits of detection because it depends on the method used for measuring the amount of dilution and also because it may be decreasing over time as measurement gets improved. In Table~1, we have reported pointwise estimation of the CDF and its standard deviation, for the Blackwood estimator and for the newly proposed estimator. We observe that at some points, the two estimators are slightly different, these points corresponding to values with both censored and uncensored measurements. As we can observe, the two estimators are the same on the right part and differ from point $15$, which is the largest value corresponding to both an exact measurement and a LOD. Below this point, the newly proposed estimator is slightly lower than the Blackwood estimator. This means that the estimate of the mean concentration will be less than the one obtained with the Blackwood estimator. As $15$ is the largest value corresponding to a LOD, the two estimators of the variance are equal for the same reason as stated above. Below this point, the estimate of the variance of the newly proposed estimator is slightly larger than the estimate of the variance of the Blackwood estimator, except for $t=x^\star_{(2)}$.
\begin{table}[h]
	\begin{center}
		\begin{tabular}{| c | c | c | c | c |}
			\hline
			$t$ & $\widehat{F}(t)$ & $ \widehat{F}^{(1)}(t)$ &$\widehat{\sigma}(\widehat{F}(t))$ & $\widehat{\sigma}(\widehat{F}^{(1)}(t))$ \\[1.5ex] \hline
			1 & 0.2981959 & 0.2799105 & 0.07438262 & 0.07541081 \\
			2 & 0.4066308 & 0.4043151 & 0.07924497  & 0.07922304  \\
			3 & 0.6235005 & 0.6199498 & 0.07582786  & 0.07644654  \\
			4 & 0.7590441 & 0.7547215 & 0.06362657  & 0.06510580  \\
			5 & 0.7820455 & 0.7816759 & 0.06125617  & 0.06159916  \\ 
			6 & 0.8280481 & 0.8276568 & 0.05555525  & 0.05598826  \\
			8 & 0.8510495 & 0.8506473 & 0.05211982 & 0.05261188 \\ 
			9 & 0.8970522 & 0.8966282 & 0.04362071  & 0.04428404  \\ 
			12 & 0.9179138 & 0.9174800 & 0.03933148  & 0.03953237  \\
			14 & 0.9387755 & 0.9383319 & 0.03424881  & 0.03449597  \\
			15 & 0.9591837 &  0.9591837 & 0.02826635  & 0.02826635  \\
			17 & 0.9795918 & 0.9795918 & 0.02019884 & 0.02019884 \\\hline
		\end{tabular}
		\caption{Pointwise estimation of the CDF and its standard deviation, for the Blackwood estimator and for the newly proposed estimator.}
	\end{center}
\end{table}

 \bibliographystyle{elsarticle-harv} 
 \bibliography{references}

\begin{thebibliography}{16}
\expandafter\ifx\csname natexlab\endcsname\relax\def\natexlab#1{#1}\fi
\providecommand{\url}[1]{\texttt{#1}}
\providecommand{\href}[2]{#2}
\providecommand{\path}[1]{#1}
\providecommand{\DOIprefix}{doi:}
\providecommand{\ArXivprefix}{arXiv:}
\providecommand{\URLprefix}{URL: }
\providecommand{\Pubmedprefix}{pmid:}
\providecommand{\doi}[1]{\href{http://dx.doi.org/#1}{\path{#1}}}
\providecommand{\Pubmed}[1]{\href{pmid:#1}{\path{#1}}}
\providecommand{\bibinfo}[2]{#2}
\ifx\xfnm\relax \def\xfnm[#1]{\unskip,\space#1}\fi
\bibitem[{Asha et~al.(2016)Asha, Elbatal and Rejeesh}]{AER}
\bibinfo{author}{Asha, G.}, \bibinfo{author}{Elbatal, I.},
  \bibinfo{author}{Rejeesh, C.J.}, \bibinfo{year}{2016}.
\newblock \bibinfo{title}{Further results on discrete mean past lifetime}.
\newblock \bibinfo{journal}{Communications in Statistics -- Theory and Methods}
  \bibinfo{volume}{45}, \bibinfo{pages}{1081--1098}.
\bibitem[{Blackwood(1991)}]{Blackwood}
\bibinfo{author}{Blackwood, L.G.}, \bibinfo{year}{1991}.
\newblock \bibinfo{title}{Analyzing censored environmental data using survival
  analysis: Single sample techniques}.
\newblock \bibinfo{journal}{Environmental Monitoring and Assessment}
  \bibinfo{volume}{18}, \bibinfo{pages}{25--40}.
\bibitem[{Gill and Johansen(1990)}]{GillJohansen90}
\bibinfo{author}{Gill, R.D.}, \bibinfo{author}{Johansen, S.},
  \bibinfo{year}{1990}.
\newblock \bibinfo{title}{A survey of product integration with a view toward
  application in survival analysis}.
\newblock \bibinfo{journal}{The Annals of Statistics} \bibinfo{volume}{18},
  \bibinfo{pages}{1501--1555}.
\bibitem[{Gillespie et~al.(2010)Gillespie, Chen, Reichert, Franzblau, Hedgeman,
  Lepkowski, Adriaens, Demond, Luksemburg and Garabrant}]{Gillespie}
\bibinfo{author}{Gillespie, B.W.}, \bibinfo{author}{Chen, Q.},
  \bibinfo{author}{Reichert, H.}, \bibinfo{author}{Franzblau, A.},
  \bibinfo{author}{Hedgeman, E.}, \bibinfo{author}{Lepkowski, J.},
  \bibinfo{author}{Adriaens, P.}, \bibinfo{author}{Demond, A.},
  \bibinfo{author}{Luksemburg, W.}, \bibinfo{author}{Garabrant, D.H.},
  \bibinfo{year}{2010}.
\newblock \bibinfo{title}{Estimating population distributions when some data
  are below a limit of detection by using a reverse {K}aplan-{M}eier
  estimator}.
\newblock \bibinfo{journal}{Epidemiology} \bibinfo{volume}{21},
  \bibinfo{pages}{S64--70}.
\bibitem[{Gomez et~al.(1992)Gomez, Juli{\`a}, Utzet and Moeschberger}]{Gomez92}
\bibinfo{author}{Gomez, G.}, \bibinfo{author}{Juli{\`a}, O.},
  \bibinfo{author}{Utzet, F.}, \bibinfo{author}{Moeschberger, M.L.},
  \bibinfo{year}{1992}.
\newblock \bibinfo{title}{Survival Analysis For Left Censored Data}.
\newblock \bibinfo{publisher}{Springer Netherlands},
  \bibinfo{address}{Dordrecht}.
\bibitem[{Gómez et~al.(1994)Gómez, Julià and Utzet}]{Gomez94}
\bibinfo{author}{Gómez, G.}, \bibinfo{author}{Julià, O.},
  \bibinfo{author}{Utzet, F.}, \bibinfo{year}{1994}.
\newblock \bibinfo{title}{Asymptotic properties of the left {K}aplan-{M}eier
  estimator}.
\newblock \bibinfo{journal}{Communications in Statistics - Theory and Methods}
  \bibinfo{volume}{23}, \bibinfo{pages}{123--135}.
\bibitem[{Helsel(2010)}]{Helsel10}
\bibinfo{author}{Helsel, D.}, \bibinfo{year}{2010}.
\newblock \bibinfo{title}{Much ado about next to nothing: incorporating
  nondetects in science}.
\newblock \bibinfo{journal}{The Annals of Occupational Hygiene}
  \bibinfo{volume}{54}, \bibinfo{pages}{257--262}.
\bibitem[{Hornung and Reed(1990)}]{HornungReed}
\bibinfo{author}{Hornung, R.W.}, \bibinfo{author}{Reed, L.D.},
  \bibinfo{year}{1990}.
\newblock \bibinfo{title}{Estimation of average concentration in the presence
  of nondetectable values}.
\newblock \bibinfo{journal}{Applied Occupational and Environmental Hygiene}
  \bibinfo{volume}{5}, \bibinfo{pages}{46--51}.
\bibitem[{Kaplan and Meier(1958)}]{KM}
\bibinfo{author}{Kaplan, E.L.}, \bibinfo{author}{Meier, P.},
  \bibinfo{year}{1958}.
\newblock \bibinfo{title}{Nonparametric estimation from incomplete
  observations}.
\newblock \bibinfo{journal}{Journal of the American Statistical Association}
  \bibinfo{volume}{53}, \bibinfo{pages}{457--481}.
\bibitem[{Millard and Deverel(1988)}]{Millard}
\bibinfo{author}{Millard, S.P.}, \bibinfo{author}{Deverel, S.J.},
  \bibinfo{year}{1988}.
\newblock \bibinfo{title}{Nonparametric statistical methods for comparing two
  sites based on data with multiple nondetect limits}.
\newblock \bibinfo{journal}{Water Resources Research} \bibinfo{volume}{24},
  \bibinfo{pages}{2087--2098}.
\bibitem[{Miller and Rupert(1981)}]{Miller}
\bibinfo{author}{Miller, J.}, \bibinfo{author}{Rupert, G.},
  \bibinfo{year}{1981}.
\newblock \bibinfo{title}{Survival analysis}.
\newblock \bibinfo{publisher}{Wiley, New York}.
\bibitem[{Pajek et~al.(2004a)Pajek, Kubala-Kukuś, Banaś, Braziewicz and
  Majewska}]{Pajek04a}
\bibinfo{author}{Pajek, M.}, \bibinfo{author}{Kubala-Kukuś, A.},
  \bibinfo{author}{Banaś, D.}, \bibinfo{author}{Braziewicz, J.},
  \bibinfo{author}{Majewska, U.}, \bibinfo{year}{2004}a.
\newblock \bibinfo{title}{Random left-censoring: a statistical approach
  accounting for detection limits in x-ray fluorescence analysis}.
\newblock \bibinfo{journal}{X-Ray Spectrometry} \bibinfo{volume}{33},
  \bibinfo{pages}{306--311}.
\bibitem[{Pajek et~al.(2004b)Pajek, Kubala-Kukuś and Braziewicz}]{Pajek04b}
\bibinfo{author}{Pajek, M.}, \bibinfo{author}{Kubala-Kukuś, A.},
  \bibinfo{author}{Braziewicz, J.}, \bibinfo{year}{2004}b.
\newblock \bibinfo{title}{Censoring: a new approach for detection limits in
  total-reflection x-ray fluorescence}.
\newblock \bibinfo{journal}{Spectrochimica Acta Part B} \bibinfo{volume}{59},
  \bibinfo{pages}{1091--1099}.
\bibitem[{Patilea and Rolin(2006)}]{PatileaRolin}
\bibinfo{author}{Patilea, V.}, \bibinfo{author}{Rolin, J.M.},
  \bibinfo{year}{2006}.
\newblock \bibinfo{title}{Product-limit estimators of the survival function
  with twice censored data}.
\newblock \bibinfo{journal}{The Annals of Statistics} \bibinfo{volume}{34},
  \bibinfo{pages}{925--938}.
\bibitem[{Popovic et~al.(2007)Popovic, Nie, Chettle and McNeill}]{Popovic}
\bibinfo{author}{Popovic, M.}, \bibinfo{author}{Nie, H.},
  \bibinfo{author}{Chettle, D.R.}, \bibinfo{author}{McNeill, F.E.},
  \bibinfo{year}{2007}.
\newblock \bibinfo{title}{Random left censoring: a second look at bone lead
  concentration measurements}.
\newblock \bibinfo{journal}{Physics in Medicine \& Biology}
  \bibinfo{volume}{52}, \bibinfo{pages}{5369}.
\bibitem[{Tressou(2006)}]{Tressou}
\bibinfo{author}{Tressou, J.}, \bibinfo{year}{2006}.
\newblock \bibinfo{title}{Nonparametric modeling of the left censorship of
  analytical data in food risk assessment}.
\newblock \bibinfo{journal}{Journal of the American Statistical Association}
  \bibinfo{volume}{101}, \bibinfo{pages}{1377--1386}.

\end{thebibliography}





\end{document}


\begin{frontmatter}



\title{Supplementary material of "New non-parametric estimators of the cumulative distribution function under time- and random-censoring"}


\author[inst1]{N. Balakrishnan}
\affiliation[inst1]{organization={Department of Mathematics and Statistics, McMaster University},
            city={Hamilton, ON},
            country={Canada, bala@mcmaster.ca} }

\author[inst2]{Ch. Paroissin}
\author[inst3]{M. Pereda Vivo}
\affiliation[inst2]{organization={Universite de Pau et des Pays de l’Adour, E2S UPPA, CNRS, LMAP},
            country={Pau, France, christian.paroissin@univ-pau.fr}}
        \affiliation[inst3]{organization={Universite de Pau et des Pays de l’Adour, E2S UPPA, CNRS, LMAP},
        	country={Pau, France, mpvivo@univ-pau.fr}}

\end{frontmatter}


{\em This project has received funding from the European Union’s Horizon 2020 research and innovation programme under the Marie Skłodowska-Curie grant agreement No 945416.}

\section{Empirical evaluation}

We conduct here a numerical comparison of the estimators, based on Monte Carlo simulations. 
We consider two cases, here. The first one deals with multiple deterministic limit of detection (LOD) values corresponding to time left-censoring scheme. The second one corresponds to random left-censoring scheme in which censored values are assumed to be realizations of random variables.

\subsection{Time left-censoring}

For the exact measurements $T_1, \ldots, T_n$, we have used the log-normal distribution with parameters $(\mu,\sigma)$, and several values of the parameters are then considered. For the LOD $C$, three possible levels, namely $0.5$, $1$, and $2$, are used in order to consider the case of multiple LODs. For each unit in the sample, one of these LODs have been selected randomly with equal probability. The sample size $n$ has been fixed to be $50$.

First, we fixed $\mu$ and then considered different values of $\sigma$. For each set of parameters, we computed the estimator in Equation~(\ref{eqn:CDF.BL.bis}) of the CDF of $T$ and the one in Equation~(\ref{estim_bis}). In order to compare the performance of these two estimators, we computed the Kolmogorov-Smirnov distance between each estimator and the exact CDF as follows :
$$
d_{ks} = \sup_{t\in\{x_{(1)}^\star,\cdots,x_{(l)}^\star\}} |F(t)-\widehat{F}(t)| 
$$
and 
$$
d_{ks}^{(1)} = \sup_{t\in\{x_{(1)}^\star,\cdots,x_{(l)}^\star\}} |F(t)-\widehat{F}^{(1)}(t)|,
$$
where $F$ is the CDF of the log-normal distribution with parameters $(\mu,\sigma)$. We repeated these steps $m = 10,000$ times and then calculated the average of the difference $d_{ks} - d_{ks}^{(1)}$ resulting from the $m$ simulations. If this average difference is positive (resp. negative), it means that globally the newly proposed estimator of $F$ is more (resp. less) accurate than the Blackwood estimator. To observe the behaviour, we first fixed $\mu$ to a given value and let $\sigma$ vary and, next, we fixed $\sigma$ to a given value and let $\mu$ vary.

\begin{center}
	\begin{figure}
		\centering
		\subfloat[]{
			\includegraphics[width=7cm]{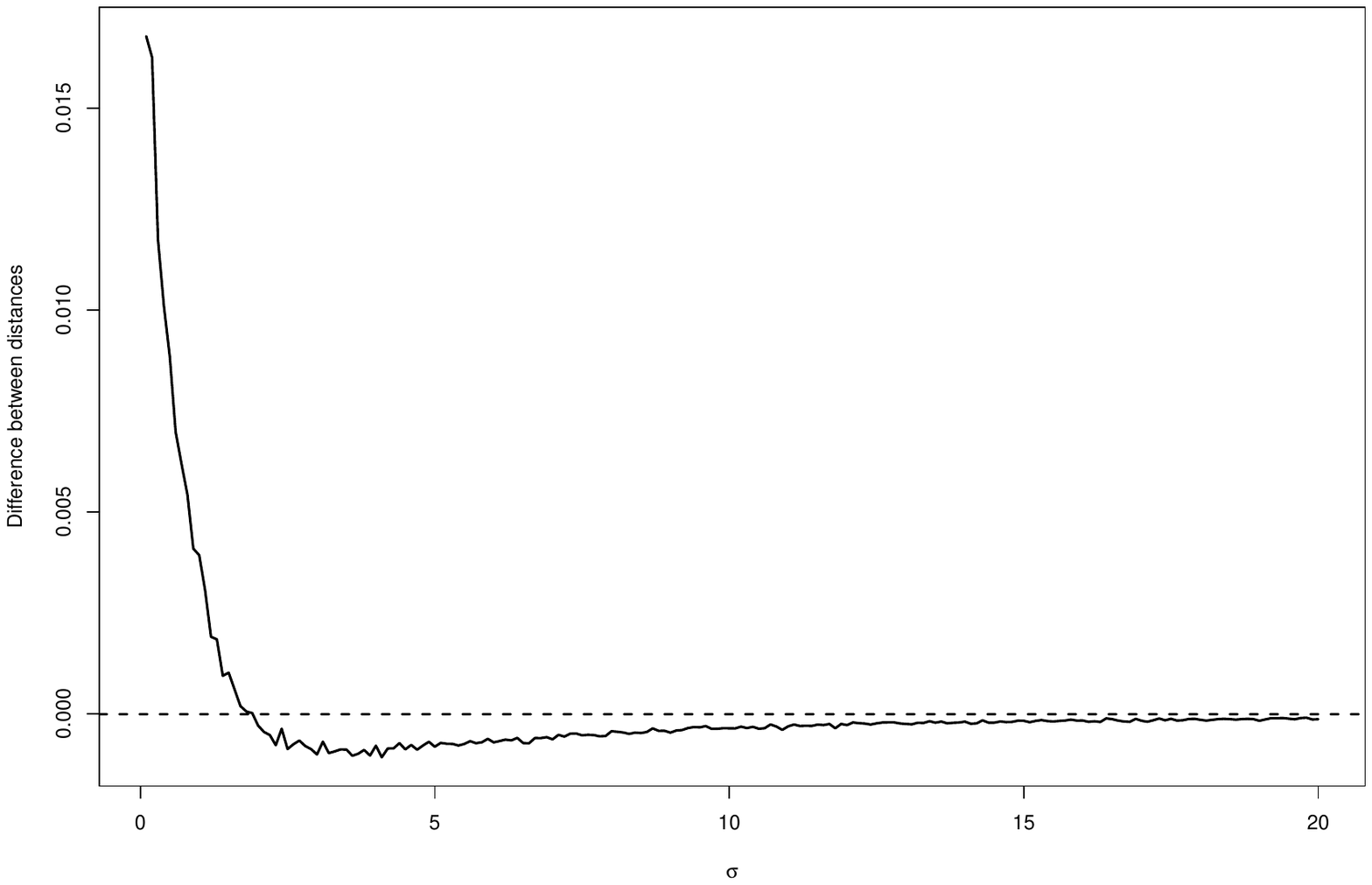}}
		\subfloat[]{
			\includegraphics[width=7cm]{type_I_mu_fixe (mean).pdf}}
		\caption{(a) Time-censoring on the left with the choices of $\mu=0$ and $\sigma\in(0,20)$; (b) Time-censoring on the left with the choices of $\mu=0$ and $\sigma\in(0,2.25)$}
		\label{fig1}
	\end{figure}
\end{center}

In Figure~\ref{fig1}, we have plotted this average difference when $\mu=0$ and $\sigma$  is varying between 0.1 and 20. We observe that the average difference is positive, next turns to be negative and then converges to zero as $\sigma$ increases.

\begin{center}
	\begin{figure}
		\centering
		\subfloat[]{
			\includegraphics[width=7cm]{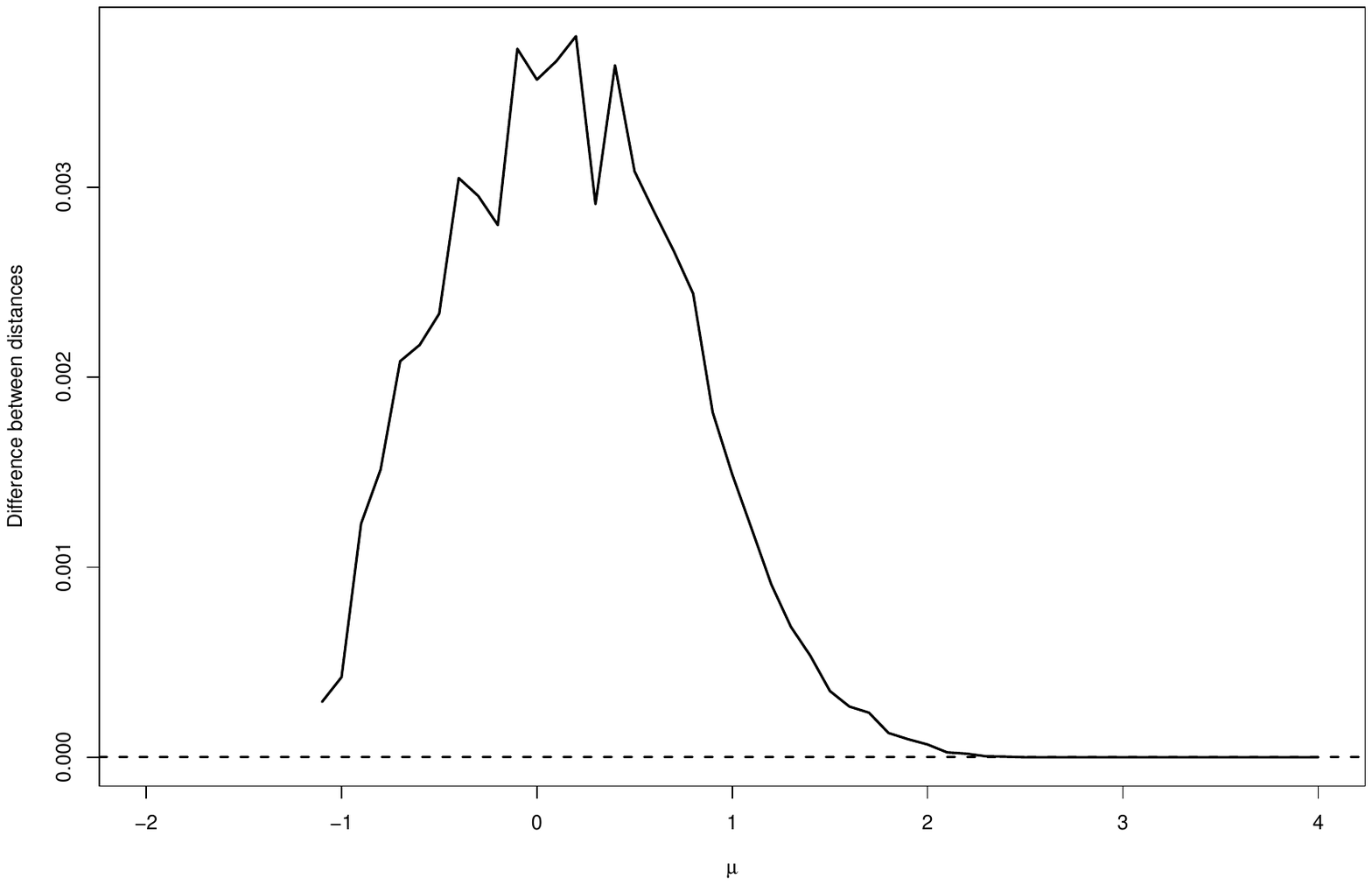}}
		\subfloat[]{
			\includegraphics[width=7cm]{type_I_sigma_fix (mean).pdf}}
		\caption{(a) Time-censoring on the left with the choices of $\sigma=1$ and $\mu\in(-2,4)$; (b) Time-censoring on the left with the choices of $\sigma=1$ and $\mu\in(-2,2)$}
		\label{fig2}
	\end{figure}
\end{center}

In Figures~\ref{fig2} and~\ref{fig3}, we have fixed $\sigma$ respectively to 1 and 2, while $\mu$ varies between -2 and 4. We observe that when $\sigma=1$ (see Figure~\ref{fig2}), the average difference  is always positive and goes to zero as $\mu$ increases. The situation is quite different when $\sigma=2$ (see Figure~\ref{fig3}). In this case, the average difference is negative, turns to be positive and then goes to zero as $\mu$ increases.

\begin{center}
	\begin{figure} 
		\centering
		\subfloat[]{
			\includegraphics[width=7cm]{type_I_sigma_fix (2).pdf}}
		\subfloat[]{
			\includegraphics[width=7cm]{type_I_sigma_fix_mean (2).pdf}}
		\caption{(a) Time-censoring on the left with the choices of $\sigma=2$ and $\mu\in(-2,4)$; (b) Time-censoring on the left with the choices of $\sigma=2$ and $\mu\in(-2,2.5)$}
		\label{fig3}
	\end{figure}
\end{center}

We thus see that when one of the two parameters of the log-normal distribution is fixed, the average difference tends to zero as the other parameter increases. This can be explained as follows. The expectation of $T$ is equal to $\exp(\mu + \sigma^2/2)$. Hence, when $\mu$ and/or $\sigma$ is increasing, the expectation of $T$ is increasing and thus one will observe more exact measurements (in other words, the probability of censoring decreases). It will make the two estimators to be closer since there will be less values under the LOD.

We obtained similar plots for $n=100$, and the sample size does not seem to affect this behaviour.

\subsection{Random left-censoring}

\begin{center}
	\begin{figure} 
		\centering
		\subfloat[]{
			\includegraphics[width=7cm]{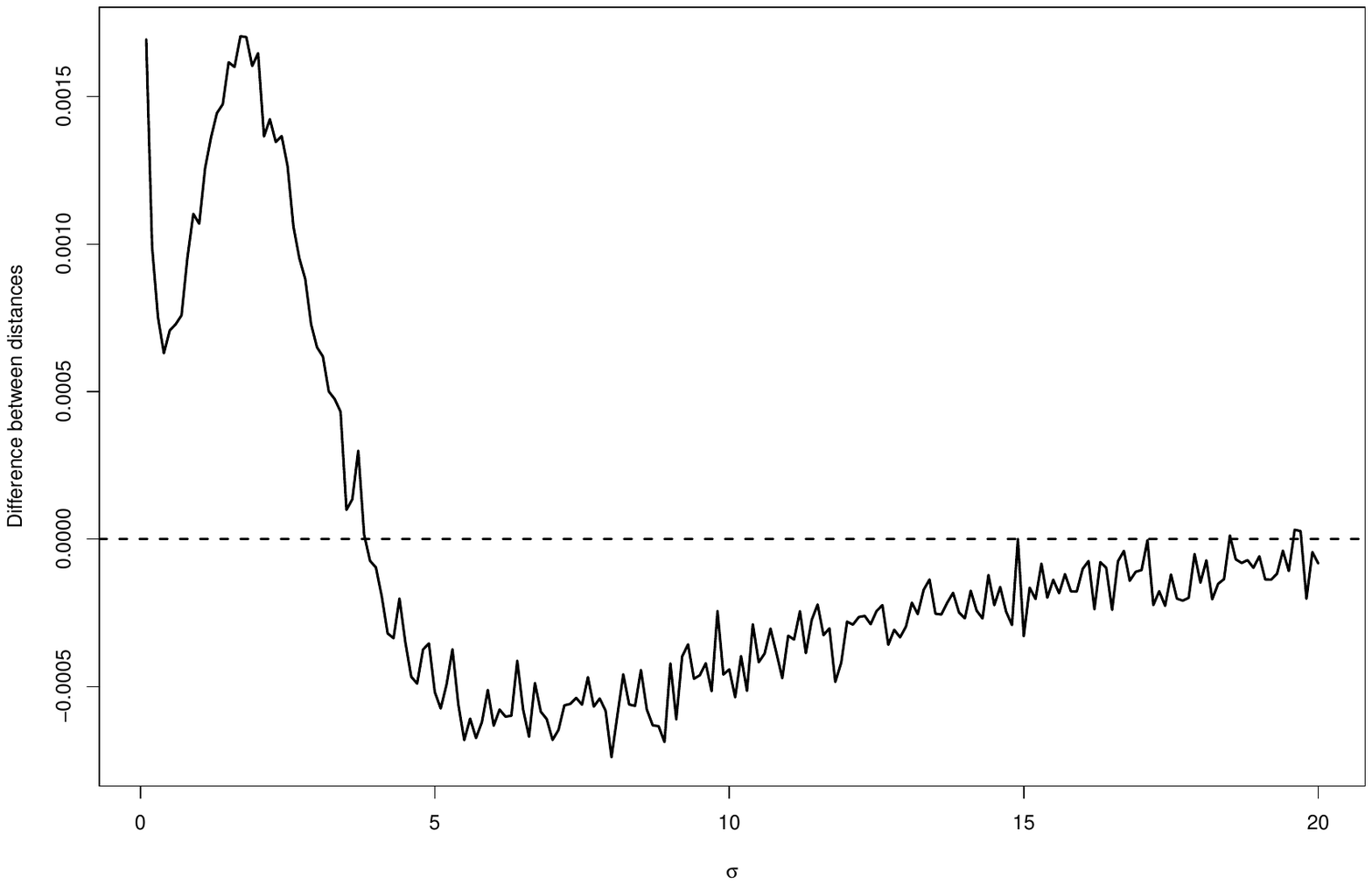}}
		\subfloat[]{
			\includegraphics[width=7cm]{type_II_mu_fixe (mean).pdf}}
		\caption{(a) Random left-censoring with the choices of $\mu=0$ and $\sigma\in(0,20)$; (b) Random left-censoring with the choices of $\mu=0$ and $\sigma\in(0,4.5)$.}
		\label{fig4A}
	\end{figure}
\end{center}

As explained above, in the case of random left-censoring, we assume that the censoring value $C$ is also random. We consider log-normal distribution with parameters $\mu_C = 0$ and $\sigma_C = 1$ for the variable $C$. As in the previous case, we will observe how the average difference evolves according to $\mu$ and $\sigma$. In Figure~\ref{fig4A}, with $\mu=0$ and $\sigma\in(0,20)$, as in the first situation of time left-censoring, we observe exactly the same behaviour (positive, negative and limit to zero).

In Figure~\ref{fig4B}, we have set $\sigma = 1$ and $\mu\in(-2,4)$. It should be noted that too small values of $\mu$ will lead to a data with only censored values in which case no estimator can be proposed. In this case, we observe the difference to be positive (still with a limit to zero).
\begin{center}
	\begin{figure} 
		\centering
		\subfloat[]{
			\includegraphics[width=7cm]{type_I(_sigma_fix.pdf}}
		\subfloat[]{
			\includegraphics[width=7cm]{type_II_sigma_fix (mean).pdf}}
		\caption{(a) Random left-censoring with the choices of $\sigma=1$ and $\mu\in(-2,4)$; (b) Random left-censoring with the choices of $\sigma=1$ and $\mu\in(-2,4)$.}
		\label{fig4B}
	\end{figure}
\end{center}





